\documentclass[12pt]{amsart}

\usepackage{amssymb,amsmath,amsthm}
\input amssym.def
\input amssym
\newtheorem{theorem}{Theorem}[section]
\newtheorem{notation}[theorem]{Notation}

\newtheorem{lemma}{Lemma}[section]

\newtheorem{remark}{Remark}[section]
\newtheorem{defi}{Definition}[section]
\newtheorem{prop}{Proposition}[section]

\newcommand{\be}{\begin{equation}}
\newcommand{\ee}{\end{equation}}

\renewcommand{\theequation}{\thesection.\arabic{equation}}
\renewcommand{\thetheorem}{\thesection.\arabic{theorem}}
\renewcommand{\theequation}{\thesection.\arabic{equation}}
\begin{document}

\title[] {Formal calculus and umbral calculus}

\author{Thomas J. Robinson}

\begin{abstract}
We use the viewpoint of the formal calculus underlying
vertex operator algebra theory to study certain aspects of the
classical umbral calculus.  We begin by calculating the
exponential generating function of the higher derivatives of a
composite function, following a very short proof which
naturally arose as a motivating computation related to a certain
crucial ``associativity'' property of an important class of vertex
operator algebras.  Very similar (somewhat forgotten) proofs had
appeared by the 19-th century, of course without any motivation
related to vertex operator algebras.  Using this formula, we derive
certain results, including especially the calculation of certain
adjoint operators, of the classical umbral calculus.  This is, roughly
speaking, a reversal of the logical development of some standard
treatments, which have obtained formulas for the higher derivatives of
a composite function, most notably Fa\`a di Bruno's formula, as a
consequence of umbral calculus.  We also show a connection between the
Virasoro algebra and the classical umbral shifts. This leads naturally
to a more general class of operators, which we introduce, and which
include the classical umbral shifts as a special case.  We prove a few
basic facts about these operators.
\end{abstract}

\maketitle

\renewcommand{\theequation}{\thesection.\arabic{equation}}
\renewcommand{\thetheorem}{\thesection.\arabic{theorem}}
\setcounter{equation}{0} \setcounter{theorem}{0}
\setcounter{section}{0}

\section{Introduction}
\setcounter{equation}{0} We present from first principles certain
aspects of the classical umbral calculus, concluding with a connection
to the Virasoro algebra.  One of our main purposes is to show
connections between the classical umbral calculus and certain central
considerations in vertex operator algebra theory.  The first major
connection is an analogue, noted in \cite{FLM}, of those authors'
original argument showing that lattice vertex operators satisfy a
certain fundamental associativity property.  Those authors observed
that this analogue amounts to a simple calculation of the higher
derivatives of a composite function, often formulated as Fa\`a di
Bruno's formula.  The philosophy of vertex operator algebra theory led
those authors to emphasize the exponential generating function of the
higher derivatives rather than the coefficients (which are easily
extracted).  That generating function was the analogue of a certain
vertex operator.  We shall show how taking this as a starting point,
one may easily (and rigorously) recover significant portions of the
classical umbral calculus of Sheffer sequences.  The main aim, part of
ongoing research, is to further develop the analogy between vertex
operator algebra theory and classical umbral calculus.  In addition, a
direct connection between the classical umbral shifts and the Virasoro
algebra (which plays a central role in vertex operator algebra theory)
is established in the second half of this paper.  Further analogies
between vertex algebra formulas and classical umbral calculus formulas
are noted in connection with this result and these motivate a
generalization of the classical umbral shifts, which we briefly
develop at the conclusion of this paper.

The classical umbral calculus has been treated rigorously in many
works following the pioneering research of Gian-Carlo Rota, such as
e.g. \cite{MR}, \cite{RKO}, \cite{Ga}, \cite{Rt}, \cite{RR},
\cite{Rm1}, \cite{Fr}, \cite{T} and \cite{Ch}.  For an extensive bibliography
through 2000 we refer the reader to \cite{BL}.  The general principle
of umbral techniques reaches far beyond the classical umbral calculus
and continues to be a subject of research (see e.g. \cite{DS},
\cite{N} and \cite{Z2}).  Our treatment involves only certain
portions of the classical umbral calculus of Sheffer sequences as
developed in \cite{Rm1}.  

There are many proofs of Fa\`a di Bruno's formula for the higher
derivatives of a composite function as well as related formulas dating
back to at least the early 19th century (see \cite{Jo} for a brief
history, as well as \cite{A}, \cite{B}, \cite{Bli}, \cite{F1},
\cite{F2}, \cite{Lu}, \cite{Me}, and \cite{Sc}).  Moreover, it is a
result that seems basic enough to be prone to showing up in numerous
unexpected places, such as in connection with vertex operator algebra
theory and also, as I recently learned from Professor Robert Wilson,
in the theory of divided power algebras, to give just one more
example.  Here, for instance, a special case of Fa\`{a} di Bruno's
formula implies that certain coefficients are combinatorial and
therefore integral, which is the point of interest since one wants a
certain construction to work over fields of finite characteristic (see
e.g. Lemma 1.3 of \cite{Wi}).  Fa\`a di Bruno's formula is purely
algebraic or combinatorial.  For a couple of combinatorial proofs we
refer the reader to \cite{Z1} and \cite{Ch}, however we shall only be
concerned with algebraic aspects of the result in this paper.

Our interest in Fa\`a di Bruno's formula is due to its appearance in
two completely separate subjects.  First, it has long well-known
connections with umbral calculus and second, perhaps more subtly, it shows up
in the theory of vertex operator algebras.  There are several umbral
style proofs of Fa\`a di Bruno's formula.  According to \cite{Jo}, an
early one of these is due to Riordan \cite{Ri1} using an argument
later completely rigorized in \cite{Rm2} and \cite{Ch}.  Perhaps even
more important, though, is the point of view taken in Section 4.1.8 of
\cite{Rm1}, where the author discusses what he calls the ``generic
associated sequence,'' which he relates to the Bell polynomials, which
themselves are closely related to Fa\`a di Bruno's formula.  The first
part of this paper may, very roughly, be regarded as showing a way to
develop some of the classical umbral calculus beginning from such
``generic'' sequences.  We also bring attention more fully to
\cite{Ch} in which the formalism of ``grammars'' and some of the
techniques quite closely resemble our approach at this stage, as I
recently became aware.

Fa\`a di Bruno's formula (in generating function form) appears in the
theory of vertex operator algebras as originally observed in
\cite{FLM}.  Briefly, Fa\`a di Bruno's formula appeared in generating
function form as an analogue, noted in \cite{FLM}, of those authors'
original argument showing that lattice vertex operators satisfy a
certain fundamental associativity property.  The work \cite{FLM} deals
with many topics, but the parts which are of interest to us have to do
with vertex operator algebra theory as well as, in particular, the
Virasoro algebra, which is a very important ingredient in vertex
operator algebra theory.  We note that although certain crucial
material from the theory of vertex operator algebras plays an
essential role in the motivation of this paper, it turns out that we
do not need explicit material directly about vertex operator algebras
for the present work. By way of the literature, we briefly mention
that the mathematical notion of vertex algebras was introduced in
\cite{B} and the variant notion of vertex operator algebra was
introduced in \cite{FLM}.  An axiomatic treatment of vertex
operator algebras was given in \cite{FHL} and a more recent
treatment was presented in \cite{LL}.  The interested reader may
consult \cite{L2} for an exposition of the history of the area.

This work began, unexpectedly, with certain considerations of the
formal calculus developed to handle some of the algebraic, and
ultimately, analytic aspects of vertex operator algebra theory.
Those considerations were related to elementary results in the
logarithmic formal calculus as developed in \cite{Mi} and \cite{HLZ}.
However, we shall not discuss the connection to the logarithmic formal
calculus here (for this see \cite{R1} and \cite{R2}) since another
more classical result stemming from vertex algebra theory turns out to
be more central to this material, namely that calculation which
amounted to a calculation of the higher derivatives of a composite
function, which was mentioned above.  For the details of this
calculation, see the introduction to Chapter 8 as well as Sections 8.3
and 8.4 of \cite{FLM} and in particular Proposition 8.3.4, formula
(8.4.32) and the comment following it.

The Virasoro algebra was studied in the characteristic $0$ case in
\cite{GF} and the characteristic $p$ analogue was introduced by
R. Block in \cite{Bl}.  Over $\mathbb{C}$ it may be realized as a
central extension of the complexified Lie algebra of polynomial vector
fields on the circle, which is itself called the Witt algebra.  A
certain crucial operator representation was introduced by Virasoro in
\cite{V} with unpublished contributions made by J.H. Weis, and the
operators of this representation play a well known and essential role
in string theory and vertex operator algebra theory (cf. \cite{FLM}).
Our connection with umbral calculus is made via one of these
operators.

Since this paper is interdisciplinary, relating ideas in vertex
operator algebra theory and umbral calculus, we have made certain
choices regarding terminology and exposition in an effort to make it
more accessible to readers who are not specialists in both of these
fields.  Out of convenience we have chosen \cite{Rm1} as a reference
for standard well-known results of umbral calculus.  A well known
feature of umbral calculus is that it is amenable to many different
recastings.  For instance, as the referee has pointed out, many of the
main classical results, recovered from our point of view in Section 4,
concerning adjoint relationships also appeared in \cite{Fr}, where what
Roman \cite{Rm1} refers to as ``adjoints'' are very nicely handled by
a certain type of ``transform.''  The change in point of view, among
other things, gives a very interesting alternative perspective on the
results and we encourage the interested reader to compare the
treatments.  However, in the interests of space, when we wish to show
the equivalence of certain of our results with the literature we will
restrict ourselves to using the notation and framework in \cite{Rm1}.

In this paper we attempt to avoid specialized vocabulary as
much as possible, although we shall try to indicate in remarks at
least some of the important vocabulary from classical works.  We shall
use the name umbral calculus or classical umbral calculus since this
seems to enjoy widespread name-recognition, but as the referee pointed
out ``finite operator calculus'' might be a more appropriate name for
much of the material such as the method in Proposition \ref{prop:FDBU}
and relevant material beginning in Section 4 of this work.  We have
also attempted to keep specialized notation to a minimum.  However,
because the notation which seems natural to begin with differs from
that used in \cite{Rm1} we do include calculations bridging the
notational gap in Section 4 for the convenience of the reader.  We
note that the proofs of the results in Section 4 are much more
roundabout than necessary if indeed those results in and of themselves
were what was sought.  The point is to show that from natural
considerations based on the generating function of the higher
derivatives of a composite function, one does indeed recover
certain results of classical umbral calculus.
  
We shall now outline the present work section-by-section.  In Section
\ref{sec:pre}, along with some basic preliminary material, we begin by
presenting a special case of the concise calculation of the
exponential generating function of the higher derivatives of a
composite function which appeared in the proof of Proposition 8.3.4 in
\cite{FLM}.  Using this as our starting point, in Section
\ref{sec:arotpafd} we then abstract this calculation and use the
resulting abstract version to derive various results of the classical
umbral calculus related to what Roman \cite{Rm1} called associated
Sheffer sequences.  The umbral results we derive in this section
essentially calculate certain adjoint operators, though in a somewhat
disguised form.  In Section \ref{sec:umbconn}, we then translate these
``disguised'' results into more familiar language using essentially
the formalism of \cite{Rm1}.  We shall also note in this section how
umbral shifts are defined as those operators satisfying what may be
regarded as an umbral analogue of the $L(-1)$-bracket-derivative
property (cf. formula (8.7.30) in \cite{FLM}).  The observation that
such analogues might be playing a role was suggested by Professor
James Lepowsky after looking at a preliminary version of this paper.

In Section \ref{sec:umbsr} we make an observation about umbral shifts
which will be useful in the last phase of the paper.

In Section \ref{sec:vir} we begin the final phase of this paper, in
which we relate the classical umbral calculus to the Virasoro algebra
of central charge 1.  Here we recall the definition of the Virasoro
algebra along with one special case of a standard ``quadratic''
representation; cf. Section 1.9 of \cite{FLM} for an exposition of
this well-known quadratic representation.  We then show how an
operator which was central to our development of the classical umbral
calculus is precisely the $L(-1)$ operator of this particular
representation of the Virasoro algebra of central charge 1.  Using a
result which we obtain in Section \ref{sec:umbsr}, we show a
relationship between the classical umbral shifts and the operator now
identified as $L(-1)$ and we then introduce those operators which in a
parallel sense correspond to $L(n)$ for $n \geq 0$.  (Strictly
speaking, by focusing on only those operators $L(n)$ with $n \geq -1$,
which themselves span a Lie algebra, the full Virasoro algebra along
with its central extension remain effectively invisible.)  We conclude
by showing a couple of characterizations of these new operators in
parallel to characterizations we already had of the umbral shifts.  In
particular we also note how the second of these characterizations,
formulated as Proposition \ref{prop:genumbshiftan}, may be regarded as
an umbral analogue of (8.7.37) in \cite{FLM}, extending an analogue
already noted concerning the $L(-1)$-bracket-derivative property.

We note also that Bernoulli polynomials have long had connections to
umbral calculus (see e.g. \cite{Mel}) and have recently appeared in
vertex algebra theory (see e.g. \cite{L1} and \cite{DLM}).  It might be
interesting to investigate further connections between the two
subjects that involve Bernoulli polynomials explicitly.

This paper is an abbreviated version of part of \cite{R3} (cf. also
\cite{R4}).  The additional material in the longer versions is largely
expository, for the convenience of readers who are not specialists.

I wish to thank my advisor, Professor James Lepowsky, as well as the
attendees (regular and irregular) of the Lie Groups/Quantum
Mathematics Seminar at Rutgers University for all of their helpful
comments concerning certain portions of the material which I presented
to them there.  I also want to thank Professors Louis Shapiro, Robert
Wilson and Doron Zeilberger for their useful remarks.  Additionally, I
would like to thank the referee for many helpful comments.

Finally, I am grateful for partial support from NSF grant
PHY0901237.
\section{Preliminaries}
\label{sec:pre}
\setcounter{equation}{0} We set up some notation and recall some
well-known and easy preliminary propositions in this section.  For a
more complete treatment, we refer the reader to the first three
sections of Chapter 8 of \cite{FLM} (cf. Chapter 2 of \cite{LL}),
while noting that in this paper we shall not need any of the material
on ``expansions of zero,'' the heart of the formal calculus treated in
those works.

We shall write $t,u,v,w,x,y,z,x_{n}, y_{m}, z_{n}$ for commuting
formal variables, where $n \geq 0$ and $m \in \mathbb{Z}$.  All vector
spaces will be over $\mathbb{C}$.  Let $V$ be a vector space.  We use
the following:
\begin{align*}
\mathbb{C}[[x]]=\biggl\{ \sum_{n \geq 0}c_{n}x^{n}|c_{n} \in
\mathbb{C} \biggr\}
\end{align*}
(formal power series), and
\begin{align*}
\mathbb{C}[x]=\biggl\{ \sum_{n \geq 0}c_{n}x^{n}|c_{n} \in \mathbb{C},
c_{n}=0 \text{ for all but finitely many } n \biggr\}
\end{align*}
(formal polynomials). 

We denote by $\frac{d}{dx}$ the formal derivative acting on either
$\mathbb{C}[x]$ or $\mathbb{C}[[x]]$.  Further, we shall frequently
use the notation $e^{\square}$ to refer to the formal exponential
expansion, where $\square$ is any formal object for which such
expansion makes sense.  By ``makes sense'' we mean that the
coefficients of the monomials of the expansion are finite objects.
For instance, we have the linear operator
$e^{w\frac{d}{dx}}:\mathbb{C}[[x,x^{-1}]] \rightarrow
\mathbb{C}[[x,x^{-1}]][[w]]$:
\begin{align*}
e^{w\frac{d}{dx}}=\sum_{n \geq
0}\frac{w^{n}}{n!}\left(\frac{d}{dx}\right)^{n}.
\end{align*}
We recall that a linear map $D$ on an algebra $A$ which satisfies
\[\begin{array}{lr}
D(ab)=(Da)b+a(Db) & \text{for all} \quad a,b \in A
\end{array}\]
is called a {\it derivation}.  Of course, the linear operator
$\frac{d}{dx}$ when acting on either $\mathbb{C}[x]$ or
$\mathbb{C}[[x]]$ is an example of a derivation.

It is a simple matter to verify, by induction for instance, the
following version of the elementary binomial theorem.  Let $A$ be an
algebra with derivation $D$.  Then for all $a$, $b \in A$, we have:
\begin{align}
 D^{n}ab&=\sum_{k+l=n}(k+l)!\frac{D^{k}a}{k!}\frac{D^{l}b}{l!}\\ 
 e^{wD}ab&=\left(e^{wD}a \right)\left(e^{wD}b\right). \qquad \qquad
\text{(the {\it automorphism property})}
\end{align}

Further, we separately state the following important special case of
the automorphism property. For $f(x),g(x) \in \mathbb{C}[[x]]$,
\begin{align*}
e^{w\frac{d}{dx}}f(x)g(x)=\left(e^{w\frac{d}{dx}}f(x)\right)
\left(e^{w\frac{d}{dx}}g(x)\right).
\end{align*}

The automorphism property shows, among other things, how the operator
$e^{w\frac{d}{dx}}$ may be regarded as a formal substitution, since,
for $n \geq 0$, we have:
\begin{align*}
e^{w\frac{d}{dx}}x^{n}=\left(e^{w\frac{d}{dx}}x\right)^{n}=(x+w)^{n}.
\end{align*}
Therefore, by linearity, we get the following polynomial formal Taylor formula.
For $p(x) \in \mathbb{C}[x]$,
\begin{align*}
e^{w\frac{d}{dx}}p(x)=p(x+w).
\end{align*}

Since the total degree of every term in $(x+w)^{n}$ is $n$, we see
that $e^{w\frac{d}{dx}}$ preserves total degree.  By equating terms
with the same total degree we can therefore extend the previous
proposition to get the following.  For $f(x) \in \mathbb{C}[[x]]$,
\begin{align}
\label{translationoperator}
e^{w\frac{d}{dx}}f(x)=f(x+w).
\end{align}

\begin{remark} \rm
We note that the identity (\ref{translationoperator}) can be derived
immediately by direct expansion as the reader may easily check.
However, in the formal calculus used in vertex operator algebra theory
it is often better to think of this minor result within a context like
that provided above.  For instance, it is often useful to regard such
formal Taylor theorems concerning formal translation operators as
representations of the automorphism property (see \cite{R1}, \cite{R2}
and Remark \ref{rem:Dan0}).
\end{remark}

We have calculated the higher derivatives of a product of two
polynomials using the automorphism property.  We next reproduce (in a
very special case, for the derivation $\frac{d}{dx}$), the quick
argument given in Proposition 8.3.4 of \cite{FLM} to calculate the
higher derivatives of the composition of two formal power series.  Let
$f(x), g(x) \in \mathbb{C}[[x]]$.  We further require that $g(x)$ have
zero constant term, so that, for instance, the composition $f(g(x))$
is always well defined.  We shall approach the problem by calculating
the exponential generating function of the higher derivatives of
$f(g(x))$.  We get
\begin{align}
\label{eq:faacalc}
e^{w\frac{d}{dx}}f(g(x))&=f(g(x+w))\nonumber\\ 
&=f(g(x)+(g(x+w)-g(x)))\nonumber\\
&=\left(e^{(g(x+w)-g(x))\frac{d}{dz}}f(z)\right)|_{z=g(x)} \nonumber\\ 
&=\sum_{n\geq
0}\frac{f^{(n)}(g(x))}{n!}\left(e^{w\frac{d}{dx}}g(x)-g(x)\right)^{n} 
\nonumber\\
&=\sum_{n \geq 0}\frac{f^{(n)}(g(x))}{n!}\left(\sum_{m \geq
1}\frac{g^{(m)}(x)}{m!}w^{m}\right)^{n}.
\end{align}  
While our calculation of the higher derivatives is not, strictly
speaking, complete at this stage (although all that remains is a
little work to extract the coefficients in powers of $w$), it is in
fact this formula which will be of importance to us, since, roughly
speaking, many results of the classical umbral calculus follow because
of it, and so we shall record it as a proposition.
\begin{prop}
\label{prop:Faa}
Let $f(x)$ and $g(x) \in \mathbb{C}[[x]]$. Let $g(x)$ have zero
constant term.  Then we have
\begin{align}
\label{eq:Faa}
e^{w\frac{d}{dx}}f(g(x))=\sum_{n \geq
0}\frac{f^{(n)}(g(x))}{n!}\left(\sum_{m \geq
1}\frac{g^{(m)}(x)}{m!}w^{m}\right)^{n}.
\end{align}  
\end{prop}
\begin{flushright} $\qed$ \end{flushright}

A derivation of Fa\`a di Bruno's classical formula may be found in
Section 12.3 of \cite{An}.  We shall not need the fully expanded
formula.

\begin{remark} \rm
\label{rem:Dan0}
The more general version of this calculation (based on a use of the
automorphism property instead of the formal Taylor theorem) appeared
in \cite{FLM} because it was related to a much more subtle and
elaborate argument showing that vertex operators associated to
lattices satisfied a certain associativity property (see \cite{FLM},
Sections 8.3 and 8.4 and in particular, formula (8.4.32) and the
comment following it).  The connection is due in part to the rough
resemblance between the exponential generating function of the higher
derivatives of a composite function in the special case $f(x)=e^{x}$
(see (\ref{Bell}) below) and ``half of'' a vertex operator.
\end{remark}

Noting that in (\ref{eq:faacalc}) we treated $g(x+w)-g(x)$ as one
atomic object suggests a reorganization.  Indeed by calling
$g(x+w)-g(x)=v$ and $g(x)=u$, the second, third and fourth lines of
(\ref{eq:faacalc}) become
\begin{align*}
f(u+v)=e^{v\frac{d}{du}}f(u)=\sum_{n \geq 0}\frac{f^{(n)}(u)}{n!}v^{n}.
\end{align*}
This is just the formal Taylor theorem, of course, and so we could
have begun here and then re-substituted for $u$ and $v$ to get the
result.  This, according to \cite{Jo}, is how the proof of U. Meyer
\cite{Me} runs.

It is also interesting to specialize to the case where $f(x)=e^{x}$, as
is often done, and indeed was the case which interested the authors
of \cite{FLM} and will interest us in later sections.  We have simply
\begin{align}
\label{Bell}
e^{w\frac{d}{dx}}e^{g(x)}=e^{g(x+w)}=e^{g(x)}e^{g(x+w)-g(x)}
=e^{g(x)}e^{^{\sum_{m \geq 1}\frac{g^{(m)}(x)}{m!}w^{m}}}.
\end{align}

\begin{remark} \rm
The generating function for what are called the Bell polynomials
(cf. Chapter 12.3 and in particular (12.3.6) in \cite{An}) easily
follows from (\ref{Bell}) using a sort of slightly unrigorous
old-fashioned umbral argument replacing $g^{(m)}$ with $g_{m}$.  See
the proof of Proposition \ref{prop:FDBU} for one way of handling such
arguments.  (The referee has pointed out that one may also rigorize
this argument with certain evaluations of umbral elements in the
umbral calculus whereas our argument in Proposition \ref{prop:FDBU} is
closer to the related finite operator calculus.) Of course, if we also
set $g(x)=e^{x}-1$, we get
$e^{w\frac{d}{dx}}e^{e^{x}-1}=e^{e^{x+w}-1}$ and setting $x=0$ is
easily seen to give the well-known result that $e^{e^{w}-1}$ is the
generating function of the Bell numbers, which are themselves the Bell
polynomials with all variables evaluated at $1$.
\end{remark}

For convenience we shall globally name three generic (up to the
indicated restrictions) elements of $\mathbb{C}[[t]]$:
\begin{align}
\label{genfun}
A(t)=\sum_{n \geq 0} A_{n}\frac{t^{n}}{n!}, \qquad
B(t)=\sum_{n \geq 1}B_{n}\frac{t^{n}}{n!}.
\qquad \text{and} \qquad
C(t)=\sum_{n \geq 0}C_{n}\frac{t^{n}}{n!},
\end{align}
where both $B_{1} \neq 0$ and $C_{0} \neq 0$ (and note the ranges of
summation).  We recall, and it is easy for the reader to check, that
$B(t)$ has a compositional inverse, which we denote by
$\overline{B}(t)$, and that $C(t)$ has a multiplicative inverse,
$C(t)^{-1}$.  We note further that since $\overline{B}(t)$ has zero
constant term, $B'(\overline{B}(t))$ is well defined, and we shall
denote it by $B^{*}(t)$.  In addition, $p(x)$ will always be a formal
polynomial and sometimes we shall feel free to use a different
variable such as $z$ in the argument of one of our generic series, so
that $A(z)$ is the same type of series as $A(t)$, only with the name
of the variable changed.
\begin{remark} \rm
We defined $B^{*}(t)=B'(\overline{B}(t))$.  As the referee pointed
out, it is also true that $B^{*}(t)=\frac{1}{\overline{B}'(t)}$, which
follows from the chain rule by taking the derivative of both sides of
$B(\overline{B}(t))=t$.
\end{remark}
\begin{remark} \rm
Series of the form $B(t)$ are sometimes called ``delta series'' in
umbral calculus, or finite operator calculus (cf. \cite{Rm1}).
\end{remark}
We shall also use the notation $A^{(n)}(t)$ for the derivatives of, in
this case, $A(t)$, and it will be convenient to define this for all $n
\in \mathbb{Z}$ to include anti-derivatives.  Of course, to make that
well-defined we need to choose particular integration constants and
only one choice is useful for us, as it turns out.

\begin{notation} \rm
For all $n \in \mathbb{Z}$, given a fixed sequence $A_{m} \in
\mathbb{C}$ for all $m \in \mathbb{Z}$, we shall define
\begin{align*}
A^{(n)}(t)=\sum_{ m \geq n}\frac{A_{m}t^{m-n}}{(m-n)!}.
\end{align*}
\end{notation}

\section{A restatement of the problem and further developments}
\setcounter{equation}{0}
\label{sec:arotpafd}
In the last section we considered the problem of calculating the
higher formal derivatives of a composite function of two formal power
series, $f(g(x))$, where we obtained an answer involving only
expressions of the form $f^{(n)}(g(x))$ and $g^{(m)}(x)$.  Because of
the restricted form of the answer it is convenient to translate the
result into a more abstract notation which retains only those
properties needed for arriving at Proposition \ref{prop:Faa}.  This
essential structure depends only on the observation that
$\frac{d}{dx}f^{(n)}(g(x))=f^{(n+1)}(g(x))(g^{(1)}(x))$ for $n \geq 0$
and that $\frac{d}{dx}g^{(m)}(x)=g^{(m+1)}(x)$ for $m \geq 1$.

Motivated by the above paragraph, we consider the algebra
$$\mathbb{C}[\dots,y_{-2},y_{-1},y_{0},y_{1},\dots, x_{1},x_{2},
\dots].$$  Then let $D$ be the unique derivation on
$\mathbb{C}[\dots,y_{-2},y_{-1},y_{0},y_{1},y_{2},\cdots, x_{1},x_{2},
\cdots]$ satisfying
\begin{align*}
\begin{array}{ll}
Dy_{i}=y_{i+1}x_{1}&  i \in \mathbb{Z}\\
Dx_{j}=x_{j+1} & j \geq 1.
\end{array}
\end{align*}
Then the question of calculating $e^{w\frac{d}{dx}}f(g(x))$ as in the
last section is seen to be essentially equivalent to calculating
\begin{align*}
e^{wD}y_{0},
\end{align*}
where we ``secretly'' identify $D$ with $\frac{d}{dx}$,
$f^{(n)}(g(x))$ with $y_{n}$ and $g^{(m)}(x)$ with $x_{m}$. We shall
make this identification rigorous in the proof of the following
proposition, while noting that the statement of said following
proposition is already (unrigorously) clear, by the ``secret''
identification in conjunction with Proposition \ref{prop:Faa}.

\begin{prop}
\label{prop:FDBU}
We have
\begin{align}
e^{wD}y_{0}&=\sum_{n \geq 0}\frac{y_{n} \left(\sum_{m \geq
1}\frac{w^{m}x_{m}}{m!}\right)^{n}}{n!} \label{eq:FDBU}.
\end{align}
\end{prop}
\begin{proof}
Let $f(x),g(x) \in \mathbb{C}[[x]]$ such that $g(x)$ has zero constant
term as in Proposition \ref{prop:Faa}.  Consider the unique algebra
homomorphism
\begin{align*}
\phi_{f,g}:\mathbb{C}[\dots,y_{-2},y_{-1},y_{0},y_{1},\dots,
x_{1},x_{2}, \dots] \rightarrow \mathbb{C}[[x]]
\end{align*}
satisfying
\[\begin{array}{lr}
\phi_{f,g}y_{i}=f^{(i)}(g(x)) &  i \in \mathbb{Z} \quad \text{and}\\ 
\phi_{f,g}x_{i}=g^{(i)}(x)  &  i \geq 1. \qquad \,\,
\end{array}\]
Then we claim that we have
\begin{align*}
\phi_{f,g} \circ D=\frac{d}{dx} \circ \phi_{f,g}.
\end{align*}
Since $\phi_{f,g}$ is a homomorphism and $D$ is a
derivation, it is clear that we need only check that these operators
agree when acting on $y_{i}$ for $i \in \mathbb{Z}$ and $x_{j}$ for $j
\geq 1$.  We get
\begin{align*}
\left(\phi_{f,g} \circ D\right)y_{i}=\phi_{f,g}(y_{i+1}x_{1})
=f^{(i+1)}(g(x))g'(x)
=\frac{d}{dx}f^{(i)}(g(x))
=\left(\frac{d}{dx} \circ \phi_{f,g}\right)y_{i}
\end{align*}
and
\begin{align*}
\left(\phi_{f,g} \circ D\right)x_{i}=\phi_{f,g}x_{i+1}
=g^{(i+1)}(x)
=\frac{d}{dx}g^{(i)}(x)
=\left(\frac{d}{dx} \circ \phi_{f,g}\right)x_{i},
\end{align*}
which gives us the claim.  Then, using the obvious extension of
$\phi_{f,g}$, by (\ref{eq:Faa}) we have
\begin{align}
\label{eq:Faa2}
\phi_{f,g}\left(e^{wD}y_{0}\right)&
=e^{w\frac{d}{dx}}\phi_{f,g}y_{0}
=e^{w\frac{d}{dx}}f(g(x))
= \phi_{f,g}\left(\sum_{n
\geq 0}\frac{y_{n}\left(\sum_{m \geq
1}\frac{w^{m}x_{m}}{m!}\right)^{n}}{n!}\right),
\end{align}
for all $f(x)$ and $g(x)$.  

Next take the formal limit as $x \rightarrow 0$ of the first and last
terms of (\ref{eq:Faa2}).  These identities clearly show that we get
identities when we substitute $f^{(n)}(0)$ for $y_{n}$ and
$g^{(n)}(0)$ for $x_{n}$ in (\ref{eq:FDBU}).  But $f^{(n)}(0)$ and
$g^{(n)}(0)$ are arbitrary and since (\ref{eq:FDBU}) amounts to a
sequence of multinomial polynomial identities when equating the
coefficients of $w^{n}$, we are done.
\end{proof}

We observe that it would have been convenient in the previous proof if
the maps $\phi_{f,g}$ had been invertible.  We provide a second proof
of Proposition \ref{prop:FDBU} using such a set-up.  This proof is
closely based on a proof appearing in \cite{Ch}.  We hope the reader
won't mind a little repetition.

\begin{proof}
(second proof of Proposition \ref{prop:FDBU})

Let $F(x)=\sum_{n \geq 0}\frac{y_{n}x^{n}}{n!}$ and $G(x)=\sum_{n \geq
1}\frac{x_{n}x^{n}}{n!}$.  Consider the unique algebra homomorphism
\begin{align*}
\psi:\mathbb{C}[\dots ,y_{-2},y_{-1},y_{0},y_{1},\dots, x_{1},x_{2},
\dots] \rightarrow \mathbb{C}[\dots ,y_{-2},y_{-1},y_{0},y_{1},\dots,
x_{1},x_{2}, \dots][[x]] 
\end{align*}
satisfying
\[\begin{array}{lr}
\psi (y_{i})=F^{(i)}(G(x)) &  i \in \mathbb{Z} \quad \text{and}\\ 
\psi (x_{i})=G^{(i)}(x)  &  i \geq 1. \qquad \,\,
\end{array}\]
Then we claim that we have
\begin{align*}
\psi \circ D=\frac{d}{dx} \circ \psi.
\end{align*}
Since $\psi$ is a homomorphism and $D$ is a derivation, it
is clear that we need only check that these operators agree when
acting on $y_{i}$ for $i \in \mathbb{Z}$ and $x_{j}$ for $j \geq 1$.
We get
\begin{align*}
\left(\psi \circ D\right)y_{i}=\psi (y_{i+1}x_{1})
=F^{(i+1)}(G(x))G'(x)
=\frac{d}{dx}F^{(i)}(G(x))
=\left(\frac{d}{dx} \circ \psi \right)y_{i}
\end{align*}
and
\begin{align*}
\left(\psi \circ D\right)x_{i}=\psi (x_{i+1})
=G^{(i+1)}(x)
=\frac{d}{dx}G^{(i)}(x)
=\left(\frac{d}{dx} \circ \psi \right)x_{i},
\end{align*}
which gives us the claim.  Then, using the obvious extension of
$\psi$,  we have
\begin{align}
\label{eq:Faa3}
\psi \left(e^{wD}y_{0}\right)&
=e^{w\frac{d}{dx}}\psi (y_{0})
=e^{w\frac{d}{dx}}F(G(x)).
\end{align}
But now we get to note that $\psi$ has a left inverse, namely setting
$x=0$, because $F^{(i)}(G(0))=y_{i}$ for $i \in \mathbb{Z}$ and
$G^{(i)}(0)=x_{i}$ for $i \geq 1$.  Thus we get
\begin{align}
\label{eq:Faa4}
e^{wD}y_{0}=\left(e^{w\frac{d}{dx}}F(G(x))\right)|_{x=0}
=F(G(x+w))|_{x=0}
=F(G(w)),
\end{align}  
which is exactly what we want.
\end{proof}

We note that our second proof of Proposition \ref{prop:FDBU} did not
depend on Proposition \ref{prop:Faa}.  Completing a natural circle of
reasoning, by using the first proof of Proposition \ref{prop:FDBU},
before invoking Proposition \ref{prop:Faa}, we had from
(\ref{eq:Faa2})
\begin{align*}
\phi_{f,g}\left(e^{wD}y_{0}\right)&
=e^{w\frac{d}{dx}}\phi_{f,g}y_{0}
=e^{w\frac{d}{dx}}f(g(x)),
\end{align*}
which by (\ref{eq:Faa4}) gives
\begin{align*}
e^{w\frac{d}{dx}}f(g(x))=\phi_{f,g}(F(G(w))),
\end{align*}
which gives us back Proposition \ref{prop:Faa}.  Thus we have shown in
a natural way how Propositions \ref{prop:Faa} and \ref{prop:FDBU} are
equivalent. 

\begin{remark} \rm
One nice aspect of our second proof of Proposition \ref{prop:FDBU},
based closely on a proof in \cite{Ch}, is that its key brings to the
fore of the argument perhaps the most striking feature of the result,
which is that the exponential generating function of higher
derivatives of a composite function is itself in the form of a
composite function.  This, of course, is an old-fashioned umbral
feature.  Furthermore, it was the form of the answer, that it roughly
resembled ``half of a vertex operator,'' which was what interested the
authors of \cite{FLM}.  This feature is also central to what follows.
\end{remark} 

We may now clearly state the trick on which (from our point of view)
much of the classical umbral calculus is based.  It is clear that if
we substitute $A_{n}$ for $y_{n}$ and $B_{n}$ for $x_{n}$ in
(\ref{eq:FDBU}) then the right-hand side will become $A(B(w))$.
Actually, it will be more interesting to substitute $xB_{n}$ for
$x_{n}$.  With this as motivation, we formally define two (for
flexibility) substitution maps.  Let $\chi_{B(t)}$ and $\psi_{A(t)}$
be the algebra homomorphisms uniquely defined by the following:
\begin{align*}
\chi_{B(t)}:\mathbb{C}[\dots,y_{-1},y_{0},y_{1},\dots,x_{1},x_{2},\cdots]
\rightarrow \mathbb{C}[\dots,y_{-1},y_{0},y_{1},\dots,x]
\end{align*}
with
\begin{align*}
\chi_{B(t)}(y_{i})&=y_{i} \qquad \qquad \quad i \in \mathbb{Z}\\
\chi_{B(t)}(x_{j})&=B_{j}x \qquad \qquad j \geq 1.
\end{align*}
and
\begin{align*}
\psi_{A(t)}:\mathbb{C}[\dots,y_{-1},y_{0},y_{1},\dots,x] \rightarrow
\mathbb{C}[x]
\end{align*}
with
\begin{align*}
\psi_{A(t)}(y_{i})&=A_{i} \qquad \qquad i \in \mathbb{Z}\\
\psi_{A(t)}(x)&=x.
\end{align*}

Then we have 
\begin{align}
\psi_{A(t)} \circ \chi_{B(t)}\left(e^{wD}y_{0}\right)
=A(xB(w)). \label{eq:genriordn}
\end{align}
To keep the notation from becoming cluttered, we shall sometimes
abbreviate $A(t)$ by simply $A$ and make other similar abbreviations
when there should be no confusion.
  
We next note that it is not difficult to explicitly calculate the
action of $\psi_{A} \circ \chi_{B} \circ e^{wD}$ on
$\mathbb{C}[\dots,y_{-1},y_{0},y_{1},\dots,x_{1},x_{2},\dots]$.
Indeed it is easy to see that we have
\begin{align}
\psi_{A} \circ \chi_{B}\circ e^{wD}y_{n}&=A^{(n)}(xB(w))
\qquad n \in \mathbb{Z} \label{eq:genriord2}\\ 
\text{\rm and } \quad \psi_{A} \circ
\chi_{B} \circ e^{wD}x_{n}&=xB^{(n)}(w) 
\,\,\, \quad \qquad n \geq 1. \label{eq:genriord3}
\end{align}
These identities determine the action completely because of the
automorphism property satisfied by $e^{wD}$.

The following series of identities (one of which is
(\ref{eq:genriordn}) repeated) is immediate from what we have shown:
\begin{align}
\psi_{A} \circ \chi_{B} \circ
e^{wD}y_{1}&=A'(xB(w)) \label{eq:1}\\ 
\psi_{A'} \circ \chi_{B} \circ e^{wD}y_{0}&=A'(xB(w))
\label{eq:2}\\ 
\frac{\partial}{\partial x} \circ \psi_{A} \circ \chi_{B} \circ
e^{wD}y_{-1}&=A(xB(w))B(w) \label{eq:3}\\
\psi_{tA(t)} \circ \chi_{B} \circ
e^{wD}y_{0}&=xB(w)A(xB(w)) \label{eq:4}\\
\psi_{A} \circ \chi_{B} \circ e^{wD}y_{0}
&=A(xB(w)) \label{eq:5}\\
\psi_{A\circ B}\circ \chi_{t} \circ e^{wD}y_{0}&=A(B(xw))
  \label{eq:6}\\ 
\frac{\partial}{\partial w}\circ
\psi_{A} \circ \chi_{B} \circ e^{wD}y_{0}&=A'(x(B(w))xB'(w)
\label{eq:7}\\ 
\psi_{B^{*}(t) A'(t)} \circ \chi_{B} \circ
e^{wD}y_{0}&=B^{*}(xB(w))A'(xB(w)). \label{eq:8}
\end{align}

We can now easily get the following proposition.
\begin{prop}
\label{prop:adjnew}
We have
\begin{enumerate}
\item 
$A'(B(w))=\left(\psi_{A} \circ \chi_{B} \circ
e^{wD}y_{1}\right)|_{x=1}= \left(\psi_{A'} \circ \chi_{B}
\circ e^{wD}y_{0}\right)|_{x=1},$
\item 
$A(B(w))B(w)=\left(
\frac{\partial}{\partial x} \circ \psi_{A} \circ \chi_{B} \circ
e^{wD}y_{-1}
\right)|_{x=1}= \left(\psi_{tA(t)}
\circ \chi_{B} \circ e^{wD}y_{0}\right)|_{x=1},$
\item 
$A(B(w))=\left(\psi_{A} \circ
\chi_{B} \circ e^{wD}y_{0}\right)|_{x=1}= \left(\psi_{A\circ
B}\circ \chi_{t} \circ e^{wD}y_{0}\right)|_{x=1},$ and
\item 
$A'(B(w))B'(w)=\frac{\partial}{\partial w}\left( \left(\psi_{A}
\circ \chi_{B}\circ e^{wD}y_{0}\right)|_{x=1}\right)=
\left(\psi_{B^{*}(t) A'(t)} \circ \chi_{B} \circ
e^{wD}y_{0}\right)|_{x=1}.$
\end{enumerate}
\end{prop}
\begin{proof}
All the identities are proved by setting $x=1$ in (\ref{eq:1}),
(\ref{eq:2}), (\ref{eq:3}), (\ref{eq:4}), (\ref{eq:5}),
(\ref{eq:6}),(\ref{eq:7}) and (\ref{eq:8}) and equating the results
pairwise as follows.  Equations (\ref{eq:1}) and (\ref{eq:2}) give
(1); equations (\ref{eq:3}) and (\ref{eq:4}) give (2); equations
(\ref{eq:5}) and (\ref{eq:6}) give (3); and equations (\ref{eq:7}) and
(\ref{eq:8}) give (4).
\end{proof}

Each of the identities in Proposition \ref{prop:adjnew} turns out to
be equivalent to the fact that a certain pair of operators are
adjoints.  In order to see this, our next task will be to put the
procedure of setting $x=1$, used in Proposition \ref{prop:adjnew},
into a context of linear functionals.  We shall do this in the next
section.

\section{Umbral connection}
\setcounter{equation}{0}
\label{sec:umbconn}
We set up a bra-ket notation following \cite{Rm1} so that we may
display some of the results obtained in the formalism there presented.

\begin{notation}\rm 
Let $f(x)=\sum_{n \geq 0}f_{n}x^{n} \in \mathbb{C}[x]$.  Then we define
\begin{align*}
\langle A(v)|f(x) \rangle =\sum_{n \geq 0}f_{n}A_{n}, 
\end{align*}
where the symbol $\langle \cdot | \cdot \rangle$ is linear in each
entry.  In particular, $\langle v^{k}/k!|x^{n} \rangle=\delta_{k,n}$,
where $\delta_{k,n}$ is the Kronecker delta.
\end{notation}
So we are now viewing $A(v)$ as a linear functional on
$\mathbb{C}[x]$.  This leads us to the notion of adjoint operators,
a key notion in the umbral calculus as presented in \cite{Rm1}.  We
shall soon show how to recover certain of the same results about
adjoints from our point of view.

\begin{defi} \rm
\label{def:adj}
We say that a linear operator $\phi$ on $\mathbb{C}[x]$ and a linear
operator $\phi^{*}$ on $\mathbb{C}[[v]]$ are adjoints exactly when,
for all $A(v)$ and for all $p(x)$, the following identity is
satisfied:
\begin{align*}
\langle \phi^{*}(A(v))|p(x) \rangle = \langle A(v)| \phi (p(x)) \rangle.
\end{align*}
\end{defi}
Of course, by linearity, it is equivalent to require that the identity
in Definition \ref{def:adj} be satisfied for $p(x)$ ranging over a
basis of $\mathbb{C}[x]$.  In addition, we extend the bra-ket notation
in the obvious way to handle elements of $\mathbb{C}[x][[w]]$
``coefficient-wise.''

\begin{prop}
\label{prop:adjch}
If $\phi$ is a linear operator on $\mathbb{C}[x]$ and $\phi^{*}$ is a
linear operator on $\mathbb{C}[[v]]$ such that
\begin{align*}
\langle \phi^{*}(A(v))|e^{xB(w)} \rangle 
= \langle A(v)|\phi \left(e^{xB(w)}\right) \rangle,
\end{align*}
for all $A(v)$ and $B(w)$, then $\phi$ and $\phi^{*}$ are adjoints.
\end{prop}
\begin{proof}
Equating coefficients of $w^{n}$ gives us the adjoint equation for a
sequence of polynomials $B_{n}(x)$ of degree exactly $n$ and arbitrary
$A(v)$.  Since the degree of $B_{n}(x)$ is $n$, these polynomials form
a basis and so the result follows by linearity.
\end{proof}

\begin{remark} \rm
\label{sheffersequencepre}
The sequence of polynomials $B_{n}(x)$ which appeared in the proof of
Proposition \ref{prop:adjch} have been called ``basic sequences'' or
sequences of ``binomial type'' (see \cite{MR}).  We shall call
them ``attached umbral sequences'' (see Definition
\ref{sheffersequence} and Remark \ref{rem:sheffersequence}).
\end{remark}

The next theorem allows us to translate our ``set $x=1$'' procedure
from Proposition \ref{prop:adjnew} into the bra-ket notation.
\begin{theorem}
\label{connection}
Let $u \in \mathbb{C}[y_{0},y_{1},\cdots, x]$ be of the form
$u=\sum_{n \geq 0}u_{n}y_{n}x^{n}$ where $u_{n} \in \mathbb{C}$.  Then
we have:
\begin{align*}
\langle A(v)|\psi_{e^{t}}(u) \rangle =(\psi_{A}(u))|_{x=1}.
\end{align*}
\end{theorem}
\begin{proof}We calculate to get:
\begin{align*}
 \langle A(v)|\psi_{e^{t}}(u) \rangle
=\langle A(v)|\sum_{n \geq 0}u_{n}x^{n} \rangle =\sum_{n \geq
0}u_{n}A_{n}=(\psi_{A}(u))|_{x=1}.
\end{align*}
\end{proof}
We may now easily observe that parts (1) and (2) of Proposition
\ref{prop:adjnew} yield adjoint relationships as the following theorem
formalizes.  Since the proofs of each statement are similar to
reasoning described below (regarding parts (3) and (4) of Proposition
\ref{prop:adjnew}) and are routine, we omit them, referring the reader
to \cite{R3} and \cite{R4} for the details.
\begin{theorem}
We have
\label{adjoints1}
\begin{enumerate}
\item 
$p(x) \in \mathbb{C}[x]$, viewed as a multiplication operator on
$\mathbb{C}[x]$ and $p(\frac{d}{dv})$ are adjoint operators.
\item 
$F(\frac{d}{dx}) \in \mathbb{C}[[\frac{d}{dx}]]$ and $F(v)$
viewed as a multiplication operator on $\mathbb{C}[[v]]$ are adjoint
operators.
\end{enumerate}
\end{theorem}
\begin{flushright} $\qed$ \end{flushright}
\begin{remark} \rm
Part (1) of Theorem \ref{adjoints1} appeared as Theorem 2.1.10 in
\cite{Rm1} and Part (2) of Theorem \ref{adjoints1} appeared as Theorem
2.2.5 in \cite{Rm1}.
\end{remark}

Parts (3) and (4) of Proposition \ref{prop:adjnew} also amount
to adjoint relationships.  By (\ref{prop:FDBU})
and Theorem \ref{connection}, we have that part (3) of Proposition
\ref{prop:adjnew} is essentially equivalent to
\begin{align}
\label{eq:umbopad}
A(B(w))&=
\langle A(v)|
\psi_{e^{t}} \circ \chi_{B} \circ e^{wD}y_{0} \rangle
=\langle A(B(v))|
\psi_{e^{t}} \circ \chi_{t} \circ e^{wD}y_{0} \rangle ,
\end{align}
which in turn, by (\ref{eq:genriordn}), gives
\begin{align}
\label{eq:umbopadj}
\langle A(v)|e^{xB(w)}\rangle =
\langle A(B(v))|e^{xw} \rangle.
\end{align}

We have therefore effectively calculated the adjoint to the
substitution map $S_{B}$ which acts by $S_{B}(g(v))=g(B(v))$ for all
$g(v) \in \mathbb{C}[[v]]$.  We simply need to make a couple of
definitions.

\begin{remark} \rm
As mentioned in the introduction, some proofs in this section are
``inefficient'' if the results are desired merely in and of
themselves.  As an example of this, we may observe that equation
(\ref{eq:umbopadj}), which essentially records a classical result as
mentioned below, is obvious once one notes that $\langle A(v) | e^{xw}
\rangle =A(w).$
\end{remark}

\begin{remark} \rm
We shall be defining certain linear operators on $\mathbb{C}[x]$ by
specifying, for instance, how they act on $e^{xw}$, which, recall,
stands for the formal exponential expansion.  Of course, by this we
mean that the operator acts only on the coefficients of $w^{n}$ $n
\geq 0$.  We have already employed similar abuses of notation with the
action of $\phi_{f,g}$ in the proof of Proposition \ref{prop:FDBU} and
with the bra-ket notation as mentioned in the comment preceding
Proposition \ref{prop:adjch}.
\end{remark}

We now recall the definition of certain ``umbral operators'';
cf. Section 3.4 in \cite{Rm1} More particularly, the umbral operator
attached to a series $B(w)$ in this work is the same as the umbral
operator for $\overline{B}(w)$ in \cite{Rm1}.  

\begin{remark} \rm
\label{rem:attachedtermin}
We shall attempt to always use the word ``attached'' in this context
to indicate the slight discrepancy of notation from Roman's \cite{Rm1}
usage regarding the switch to the compositional inverse.
\end{remark}

\begin{defi} \rm
\label{umbraloperator}
We define the umbral operator attached to $B(w)$ to be the unique linear map
$\theta_{B} : \mathbb{C}[x] \rightarrow \mathbb{C}[x]$ satisfying:

\begin{align*}
\theta_{B}e^{xw}=e^{xB(w)}.
\end{align*}
\end{defi}

\begin{theorem}
\label{adjoints2}
We have that $S_{B}$ and $\theta_{B}$ are adjoint operators.
\end{theorem}
\begin{proof}
The result follows from Proposition \ref{prop:adjch} and (\ref{eq:umbopadj}).
\end{proof}

\begin{remark} \rm
\label{rem:refadj2}
Theorem \ref{adjoints2} essentially appeared as Theorem 3.4.1 in
\cite{Rm1}, although in this work we have chosen some different
characterizations of certain objects as definitions as discussed in
the introduction.  It is not difficult to tie up all the relevant
information.  For more details see \cite{R3}.
\end{remark}

By (\ref{prop:FDBU}) and Theorem \ref{connection} we have that part
(4) of Proposition \ref{prop:adjnew} is essentially equivalent to
\begin{align*}
A'(B(w))B'(w)&=\frac{\partial}{ \partial w}
\langle A(v)|\psi_{e^{t}} \circ
\chi_{B} \circ e^{wD}y_{0}\rangle 
=\langle B^{*}(v)A'(v)|\psi_{e^{t}}
\circ \chi_{B} \circ e^{wD}y_{0} \rangle,
\end{align*}
which in turn, by (\ref{eq:genriordn}), gives
\begin{align}
\label{eq:umbshadj}
\frac{\partial}{\partial w}
\langle A(v)|e^{xB(w)} \rangle 
&=\langle B^{*}(v)A'(v)|e^{xB(w)} \rangle
\quad \Leftrightarrow \nonumber\\ 
\langle A(v)|\frac{\partial}{\partial
w}e^{xB(w)} \rangle
&=\langle B^{*}(v)A'(v)|e^{xB(w)} \rangle.
\end{align}

We now recall the definition of a certain important special class of
``Sheffer shifts''; cf. Section 3.6 in \cite{Rm1}.  More particularly,
the umbral shift attached to a series $B(w)$ in this work is the
same as the umbral shift for $\overline{B}(w)$ in \cite{Rm1}.

\begin{remark} \rm
Regarding the word ``attached,'' in this context see Remark
\ref{rem:attachedtermin}.
\end{remark}

\begin{defi} \rm
\label{def:umbshift}
For each $B(w)$, let $D_{B}:\mathbb{C}[x] \rightarrow \mathbb{C}[x]$
be the unique linear map satisfying
\begin{align}
\label{Dderiv}
D_{B}e^{xB(w)}=\frac{\partial}{\partial w}e^{xB(w)}.
\end{align}
We call $D_{B}$ the {\it umbral shift} attached to $B(w)$.  
\end{defi}

\begin{remark} \rm
\label{rem:Dan1}
As discussed in the Introduction and Remark \ref{rem:Dan0}, the
authors of \cite{FLM} were concerned with the exponential generating
function of the higher derivatives of a composite function, because it
roughly resembled ``half of'' a vertex operator.  Following this
analogy, we might say that $\psi_{e^{t}} \circ \chi_{B} \circ
e^{wD}y_{0}$ is an umbral analogue of (``half of'') a vertex operator.
Having made this analogy, one can see how using (\ref{Dderiv}), we
have defined ``attached'' umbral shifts as those operators satisfying an
analogue of the $L(-1)$-bracket-derivative property, which is stated
as the equality of the first and third expressions of formula (8.7.30)
in \cite{FLM}.  See also Remark \ref{rem:Dan2}.
\end{remark}

\begin{theorem}
\label{adjoints3}
We have that $D_{B}$ and $B^{*}(v) \circ \frac{d}{dv}$ are adjoints.
\end{theorem}
\begin{proof}
The result follows from Proposition \ref{prop:adjch} and (\ref{eq:umbshadj}).
\end{proof}

\begin{remark} \rm
Theorem \ref{adjoints3} essentially appeared as part of Theorem 3.6.1
in \cite{Rm1}, where the author of that work had already, in addition,
shown that the operators $B^{*}(t) \circ \frac{d}{dt}$ are exactly the
surjective derivations on $\mathbb{C}[[t]]$, a routine matter once we
note that $B^{*}(t)$ is an arbitrary element of $\mathbb{C}[[t]]$
having a multiplicative inverse.  We also mention a similar caveat for
the reader regarding different choices of definitions between the
present work and \cite{Rm1} just as discussed in Remark
\ref{rem:refadj2} and in the Introduction.
\end{remark}

In closing this section we note obvious characterizations of the
attached umbral operators and attached umbral shifts in terms of the
coefficients of their generating function definitions.  For this it is
convenient for us to recall the definition of attached umbral
sequences; cf. Section 2.3 and Theorem 2.3.4 in particular in
\cite{Rm1}.  We note that the umbral sequence ``attached'' to a series
$B(w)$ in this work is the same as the Sheffer sequence associated to
$\overline{B}(w)$ in \cite{Rm1}.

\begin{remark} \rm
Regarding the use of the word ``attached'' in conjunction with umbral
sequence, see Remark \ref{rem:attachedtermin}.
\end{remark}

\begin{defi} \rm
\label{sheffersequence}
We define the sequence of polynomials $B_{n}(x)$, the umbral sequence
attached to $B(w)$, to be the unique sequence satisfying the following:
\begin{align*}
e^{xB(w)}=\sum_{n \geq 0}\frac{x^{n}B(w)^{n}}{n!}=\sum_{n \geq
0}\frac{B_{n}(x)w^{n}}{n!}.
\end{align*}
\end{defi}

\begin{remark} \rm
\label{rem:sheffersequence}
We recall that the attached umbral sequences already appeared
explicitly (though, of course not by name) in the proof of Proposition
\ref{prop:adjch} (see Remark \ref{sheffersequencepre}).
\end{remark}

\begin{prop}
We have that $\theta_{B}:\mathbb{C}[x] \rightarrow \mathbb{C}[x]$, the
umbral operator attached to $B(w)$, is characterized as the unique
linear map satisfying:
\begin{align*}
\theta_{B}x^{n}=B_{n}(x).
\end{align*}
\end{prop}
\begin{flushright} $\qed$ \end{flushright}

\begin{prop}
\label{prop:attsheffpolcal}
We have that $D_{B}:\mathbb{C}[x] \rightarrow
\mathbb{C}[x]$, the umbral shift attached to $B(w)$, is characterized
as the unique linear map satisfying:
\begin{align*}
D_{B}B_{n}(x)=B_{n+1}(x).
\end{align*}
\end{prop}
\begin{flushright} $\qed$ \end{flushright}

\section{Umbral shifts revisited}
\setcounter{equation}{0} 
\label{sec:umbsr}
In this section we shall show a
characterization of the attached umbral shifts which will be useful
in Section \ref{sec:umbvir2}.  We begin by (for temporary convenience)
generalizing Definition \ref{def:umbshift}.
\begin{defi} \rm
\label{def:genumbshift}
For each $A(t)$ and $B(t)$, let $D^{A}_{B}:\mathbb{C}[x]
\rightarrow \mathbb{C}[x]$ be the unique linear map satisfying
\begin{align*}
D^{A}_{B}e^{xB(t)}=\frac{\partial}{\partial t}A(xB(t)).
\end{align*}
\end{defi}

Recalling the identities (\ref{eq:genriord2}) and (\ref{eq:genriord3}),
we note that
\begin{align*}
\psi_{A} \circ \chi_{B} \circ De^{wD}y_{0}
&= \psi_{A} \circ \chi_{B} \circ e^{wD}y_{1}x_{1}\\ 
&=\psi_{A} \circ \chi_{B} \circ \left(e^{wD}y_{1}\right)
\left(e^{wD}x_{1}\right)\\ 
&=A'(xB(w))xB'(w),
\end{align*}
so that
\begin{align*}
A'(xB(w))xB'(w)
&=\psi_{A} \circ \chi_{B} \circ De^{wD}y_{0}\\
&=\frac{\partial}{\partial w} \left(\psi_{A} \circ \chi_{B} \circ
e^{wD}y_{0}\right)\\
&=D^{A}_{B} \circ \psi_{e^{t}} \circ
\chi_{B} \circ e^{wD}y_{0}.
\end{align*}
Recalling that $e^{wD}$ stands for the formal exponential Taylor
series, and extracting the coefficients in $w^{n}/n!$ for $n \geq 0$
from the second and fourth terms from the above identity yields:
\begin{align*}
D^{A}_{B} \circ \psi_{e^{t}} \circ \chi_{B} \circ D^{n}y_{0}
=\psi_{A} \circ \chi_{B} \circ D^{n+1}y_{0}.
\end{align*}
Furthermore, because $\psi_{e^{t}} \circ \chi_{B} \circ
D^{n}y_{0}$ is a polynomial of degree exactly $n$, this
formula characterizes the maps $D^{A}_{B}$.  

Although we have briefly generalized the definition for the attached
umbral shifts (in order to fit more closely with our calculations from
Section \ref{sec:arotpafd}), the previous identity shows how it is natural to
restrict our attention to the attached umbral shifts, and it is this
case that will later interest us anyway.  We may now state the
characterization of the attached umbral shifts mentioned in the
introduction to this section.

\begin{remark} \rm
While we only temporarily generalized the definition for attached
umbral shifts, as the referee has pointed out, it might be nice to
investigate whether extensions of standard umbral calculus
calculations could be developed for the operators $D^{A}_{B}$.
\end{remark}
\begin{prop}
\label{prop:umbshiftchar}
The attached umbral shift, $D_{B}:\mathbb{C}[x] \rightarrow
\mathbb{C}[x]$ is the unique linear map satisfying
\begin{align*}
D_{B} \circ \psi_{e^{t}} \circ \chi_{B} \circ
D^{n}y_{0} =\psi_{e^{t}} \circ \chi_{B} \circ
D^{n+1}y_{0},
\end{align*}
for all $n \geq 0$.
\end{prop}
\begin{flushright} $\qed$ \end{flushright}

\begin{remark} \rm
Proposition \ref{prop:umbshiftchar} was announced, together with a
more direct proof, as Proposition 6.1 in \cite{R2}.
\end{remark}

\begin{remark} \rm
The classical umbral calculus can be considered to be the study of
Sheffer sequences through ``umbral'' techniques.  So far we have only
considered a special case, attached umbral sequences.  The general
case may be obtained easily by using the results of the special case
and this is probably the shortest route given the efforts we have
already made.  Alternatively one could give a complete parallel
development.  Recall that our approach began by calculating the higher
derivatives of a formal composite function.  That is, letting $f(x),
g(x) \in \mathbb{C}[[x]]$ such that the constant term of $g(x)$ is
zero, we began calculating the higher derivatives of $f(g(x))$.  Let
$h(x) \in \mathbb{C}[[x]]$.  Then the general theory follows by the
parallel argument with the starting point of calculating the higher
derivatives of the product $h(x)f(g(x))$.  The interested reader may
see this worked out in detail in both \cite{R3} and \cite{R4}.
\end{remark}

\section{Connection with the Virasoro algebra}
\setcounter{equation}{0} 
\label{sec:vir}
Our goal in this section will be to show that an operator closely
related to the derivation $D$ (which appeared in Section
\ref{sec:arotpafd}) is one of the standard quadratic representations
of the $L(-1)$ operator of the Virasoro algebra.

Recall that we began our main investigation by calculating the higher
derivatives of the composition of two formal power series $f(x)$ and
$g(x)$, where the constant term of $g(x)$ was required to be $0$,
following a proof given in \cite{FLM}.  In fact, the case that
interested the authors in \cite{FLM} was when $f(x)=e^{x}$.  It is not
difficult to specialize our arguments to this case.  When we abstract,
we get the following set-up: Consider the vector space
$y\mathbb{C}[x_{1},x_{2},x_{3}, \dots]$ where $ x_{j}$ for $j \geq 1$
are commuting formal variables.  Then let $\mathcal{D}$ be the unique
derivation on $y\mathbb{C}[x_{1},x_{2},x_{3}, \dots]$ satisfying
\begin{align}
\label{virD}
\begin{array}{ll}
\mathcal{D}y=yx_{1}&  \\
\mathcal{D}x_{j}=x_{j+1} & j \geq 1.
\end{array}
\end{align}
The question of calculating $e^{w\frac{d}{dx}}e^{g(x)}$ is seen to be
essentially equivalent to calculating
\begin{align*}
e^{w\mathcal{D}}y,
\end{align*}
where we ``secretly'' identify $\mathcal{D}$ with $\frac{d}{dx}$,
$e^{g(x)}$ with $y$ and $g^{(m)}(x)$ with $x_{m}$.  It is clear by our
identification, and rigorously as an easy corollary of Proposition
\ref{prop:FDBU}, that:
\begin{align}
e^{w\mathcal{D}}y&=
\sum_{n \geq 0}\frac{y\left(\sum_{m \geq 1}
\frac{w^{m}x_{m}}{m!}\right)^{n}}{n!}=
ye^{^{\sum_{m \geq 1}\frac{w^{m}x_{m}}{m!}}}.
\end{align}

We note that 
\begin{align*}
\mathcal{D}=x_{1}y\frac{\partial}{\partial
y}+x_{2}\frac{\partial}{\partial x_{1}}+x_{3}\frac{\partial}{\partial
x_{2}}+\cdots\\ 
=x_{1}y\frac{\partial}{\partial y}+\sum_{k \geq
1}x_{k+1}\frac{\partial}{\partial x_{k}}.
\end{align*}

We next shall switch gears in order to recall certain basics about the
Virasoro algebra using operators arising from certain Heisenberg Lie
algebras.  We follow (a variant of) the exposition of this well-known
material in \cite{FLM}.  Let $\mathfrak{h}$ be the one-dimensional
abelian (complex) Lie algebra with basis element $h$.  We define a
nonsingular symmetric bilinear form on $\mathfrak{h}$ by $(ah,bh)=ab$
for all $a,b \in \mathbb{C}$.  We recall the (particular) affine
Heisenberg Lie algebra $\widehat{\mathfrak{h}}$ which is the vector
space
\begin{align*}
\widehat{\mathfrak{h}}=
\mathfrak{h} \otimes \mathbb{C}[t,t^{-1}] \oplus \mathbb{C}c,
\end{align*}
with Lie brackets determined by
\begin{align*}
[ah\otimes t^{m},bh \otimes
t^{n}]=(ah,bh)m\delta_{m+n,0}c=abm\delta_{m+n,0}c,
\end{align*}
where $c$ is central and $\delta$ is the Kronecker delta.  

We may realize $\widehat{\mathfrak{h}}$ as differential and
multiplication operators on a space with infinitely many variables as
follows.  We consider the space $y\mathbb{C}[x_{1},x_{2},x_{3},\dots]$
and make the following identification:
\begin{align*}
h \otimes t^{n}=
\left\{
\begin{array}{lr}
\alpha(-n)x_{-n}& n < 0\\
\beta(n)\frac{\partial}{\partial x_{n}} & n > 0\\
y \frac{\partial}{\partial y}& n=0,
\end{array}
\right.
\end{align*}
where $\alpha(n),\beta(n) \in \mathbb{C}$ for $n \geq 1$ and we
identify $c$, the central element, with the multiplication by identity
operator.  Of course, in this setting $y \frac{\partial}{\partial y}$
is a fancy name for the identity operator, but we wrote it this way so
that it appears explicitly as a derivation.  It is easy to see that
\begin{align*}
\left[\alpha(n)x_{n},\beta(n)\frac{\partial}{\partial x_{n}}\right]
=-\alpha(n)\beta(n),
\end{align*}
with all other pairs of operators commuting.  Thus our identification
gives a representation of the Heisenberg Lie algebra exactly when we
require that for all $n \geq 1$
\begin{align*}
\alpha(n)\beta(n)=n.
\end{align*}  
For this representation we shall sometimes use the notation $h(n)$ to
denote the image of $h\otimes t^{n}$.

\begin{defi} \rm
The Virasoro algebra of central charge 1 is the Lie algebra generated
by basis elements $1$, a central element, and $L(n)$ for $n \in
\mathbb{Z}$ which satisfy the following relations for all $m,n \in
\mathbb{Z}$.
\begin{align}
\label{vir}
[L(m),L(n)]=(m-n)L(m+n)+(1/12)(m^{3}-m)\delta_{m+n,0}.
\end{align}
\end{defi}

\begin{remark} \rm
The Virasoro algebra is a central extension of the Witt algebra, the
Lie algebra of the derivations of Laurent polynomials in a single
formal variable.  The correspondence can be seen by identifying $L(n)$
with $-t^{n+1}\frac{d}{dt}$.  In fact, (cf. Proposition 1.9.4 in
\cite{FLM}) the Virasoro algebra with a general central element is the
unique, up to isomorphism, one dimensional central extension of the
Witt algebra.
\end{remark}

We state the following theorem without proof.  It is a special case,
for instance of Theorem 1.9.6 in \cite{FLM} where a complete proof is
provided.
    
\begin{theorem}
\label{prop:vir}
  The operators
\begin{align}
L(n)&=\frac{1}{2}\sum_{k \in \mathbb{Z}}h(n-k)h(k) 
\qquad n \neq 0 \qquad \text{and} \label{quadvir}\\
L(0)&=\frac{1}{2}\sum_{k \in \mathbb{Z}}h(-|k|)h(|k|) \label{quadvir2}
\end{align}
give a representation of the Virasoro algebra of central charge $1$,
\end{theorem}

The space $y\mathbb{C}[x_{1},x_{2},x_{3},\dots]$ is obviously a module
for the Virasoro algebra.  It is graded by $L(0)$ eigenvalues, which
are called {\it weights}.  In the literature, such a module is often
called a lowest weight module; this module has $y$ as a lowest weight
vector.

We have
\begin{align*}
L(0)&=\frac{1}{2}h(0)^{2}+h(-1)h(1)+h(-2)h(2)+\cdots\\
&=\frac{1}{2}y\frac{\partial}{\partial y} \circ
y\frac{\partial}{\partial
y}+\alpha(1)\beta(1)x_{1}\frac{\partial}{\partial
x_{1}}+\alpha(2)\beta(2)x_{2}\frac{\partial}{\partial x_{2}}+\cdots\\
\end{align*}
so that $L(0)y=\frac{1}{2}y$.  Thus the lowest weight of the module
is $\frac{1}{2}$.

We may now show that by an appropriate (unique) choice of $\alpha(n)$
and $\beta(n)$ we get $\mathcal{D}=L(-1)$.  We have
\begin{align*}
L(-1)&=h(-1)h(0)+h(-2)h(1)+h(-3)h(2)+\dots\\
&=\alpha(1)x_{1}y\frac{\partial}{\partial y}
+\alpha(2)\beta(1)x_{2}\frac{\partial}{\partial x_{1}}
+\alpha(3)\beta(2)x_{3}\frac{\partial}{\partial x_{2}}+ \dots.
\end{align*}
Therefore it is clear that in order to have $\mathcal{D}=L(-1)$, we
need exactly that
\begin{align*}
\alpha(1)&=1 \qquad \text{and}\\
\alpha(n+1)\beta(n)&=1 \quad n \geq 1,
\end{align*}
where we recall that we already have the restriction that
$\alpha(n)\beta(n)=n$ for all $n \geq 1$.  These two sets of
restrictions imply that
\begin{align*}
(n+1)\beta(n)&=\beta(n+1)  \quad n \geq 1\\
\beta(1)&=1,
\end{align*}
so that 
\begin{align*}
\beta(n)&=n!\\
\alpha(n)&=\frac{1}{(n-1)!},
\end{align*}
for all $n \geq 1$, is the unique solution.  We record this as a
proposition.

\begin{prop}
The operator $L(-1)$, given by (\ref{quadvir}), with (and only with)
both $\alpha(n)=\frac{1}{(n-1)!}$ and $\beta(n)=n!$, is identical to
the operator $\mathcal{D}$, given by (\ref{virD}).
\end{prop}
\begin{flushright} $\qed$ \end{flushright}

For the remainder of this paper we shall assume that
$\alpha(n)=\frac{1}{(n-1)!}$ and $\beta(n)=n!$.

\section{Umbral shifts revisited and generalized}
\label{sec:umbvir2}
\setcounter{equation}{0} We shall continue to consider the space
$y\mathbb{C}[x_{1},x_{2},x_{3},\dots]$ as in the previous section, and
similarly to some of our previous work, such as in Section
\ref{sec:arotpafd}, we shall consider certain substitution maps.  Let
$\phi_{B(t)}$ denote the following algebra homomorphism.
\begin{align*}
\phi_{B}:y\mathbb{C}[x_{1},x_{2},x_{3},\dots] \rightarrow \mathbb{C}[x]
\end{align*} 
with
\begin{align*}
\phi_{B(t)}x_{j}&=B_{j}x  \qquad j \geq 1 \\
\text{and} \qquad \phi_{B(t)}y&=1.
\end{align*}
Then we have
\begin{align*}
\phi_{B} \circ e^{w\mathcal{D}}y=\phi_{B} \circ e^{wL(-1)}y=e^{xB(w)}.
\end{align*}
In light of Proposition \ref{prop:umbshiftchar}, it is routine to show
the following.
\begin{theorem}
\label{theorem:umbvir}
The attached umbral shift, $D_{B}:\mathbb{C}[x] \rightarrow
\mathbb{C}[x]$ is the unique linear map satisfying
\begin{align*}
D_{B} \circ \phi_{B} \circ L(-1)^{n}y 
=\phi_{B} \circ L(-1)^{n+1}y,
\end{align*}
for all $n \geq 0$.
\end{theorem}
\begin{flushright} $\qed$ \end{flushright}

With Theorem \ref{theorem:umbvir} as motivation, we make the following
definition.
\begin{defi} \rm
\label{def:genattshef}
For $m \geq -1$ we define the operators $\mathcal{D}_{B}(m):
\mathbb{C}[x] \rightarrow \mathbb{C}[x]$ to be the unique linear maps
satisfying
\begin{align*}
\mathcal{D}_{B}(m) \circ \phi_{B} \circ L(-1)^{n}y 
=\phi_{B} \circ L(m)L(-1)^{n}y,
\end{align*}
for all $n \geq 0$.
\end{defi}
Of course, $\mathcal{D}_{B}(-1)=D_{B}$.  These operators are
well-defined because $\phi_{B} \circ L(-1)^{n}y$ has degree exactly
$n$.  In fact, $\phi_{B} \circ L(-1)^{n}y=B_{n}(x)$, the umbral
polynomial attached to $B(t)$.  We call the operators
$\mathcal{D}_{B}(n)$ generalized attached umbral shifts.

We would like to use the Virasoro relations to help compute the
generalized attached umbral shifts.  We begin with the following lemma.
\begin{lemma}
\label{fsubmlemma}
There exist rational numbers $f_{m}(n)$ such that
\begin{align*}
L(m)L(-1)^{n}y=f_{m}(n)L(-1)^{n-m}y,
\end{align*}
for all $m \geq -1$, $n \geq 0$ such that $n \geq m$.
\end{lemma}
\begin{proof}
We will need that $L(m)L(-1)^{n}y=0$ when $m > n$.  This is really due
to the weights of the vectors, but in this paper we shall proceed, in
just this special case, with an elementary induction argument.  We
induct on $n$.  For $n=0$ this follows essentially because $y$ is a
lowest weight vector, but even without considering weights it is easy
to directly see given the definition of the operators.  By induction
(used twice) we have
\begin{align*}
L(m)L(-1)^{n}y&=L(-1)L(m)L(-1)^{n-1}y+[L(m),L(-1)]L(-1)^{n-1}y\\
&=(m+1)L(m-1)L(-1)^{n-1}y\\
&=0.
\end{align*}
We may now focus on the main argument.  We establish the boundary
cases.  Letting $m=-1$ we easily check that $f_{-1}(n)=1$.  The other
boundary is $m=n$.  We shall use another intermediate induction to
establish this case.  Our base case then is $m=n=0$ for which it is
easy to check that $f_{0}(0)=1/2$.  We also have
\begin{align*}
L(n)L(-1)^{n}y&=L(-1)L(n)L(-1)^{n-1}y+[L(n),L(-1)]L(-1)^{n-1}y\\
&=(n+1)L(n-1)L(-1)^{n-1}y,
\end{align*}
so that inducting on $n$ we get our result. Moreover we now have the
recurrence
\begin{align}
\label{diagrec}
f_{n}(n)=(n+1)f_{n-1}(n-1) \qquad n \geq 1,
\end{align}
with, as we have seen, the boundary $f_{0}(0)=1/2$.

For our main argument we induct on $m+n$.  We have already checked the
base case.  We then have by induction and using the Virasoro relations
that for the remaining cases $m \geq 0$ and $n > m$, we have
\begin{align*}
L(m)L(-1)^{n}y&=L(-1)L(m)L(-1)^{n-1}y+[L(m),L(-1)]L(-1)^{n-1}y\\
&=f_{m}(n-1)L(-1)^{n-m}y+(m+1)L(m-1)L(-1)^{n-1}y\\
&=\left(f_{m}(n-1)+(m+1)f_{m-1}(n-1)\right)L(-1)^{n-m}y.
\end{align*}
Therefore, not only do the values $f_{n}(m)$ exist but we have a
recurrence for them 
\begin{align}
\label{rec}
f_{m}(n)=f_{m}(n-1)+(m+1)f_{m-1}(n-1).
\end{align}
\end{proof}

In the last proposition we found a recurrence for certain values
$f_{m}(n)$.  We could extend the range of $m$ and $n$ to include all
$m \geq -1$ and $n \in \mathbb{Z}$ and define $f_{m}(n)$ to be the
solution to the recurrence equation found above which coincides when
$n \geq 0$ and $n > m$ with the values already defined.  In fact, it
is easy to find a simpler boundary condition yielding the desired
solution other than using the boundary with $m=n$.  It is easy to see
that instead we may specify that $f_{m}(0)=0$ for $m \geq 1$, by
considering (\ref{diagrec}), which shows that we may specify $0$'s
below the diagonal.  Further, it is easy to see from this recurrence
equation with given boundary, that $f_{m}(n)$ is an integer for $m
\neq 0$ and that $f_{0}(n)$ are half integers.  It is also easy to see
from this recurrence, by induction on $n$, that we have for $n \geq 0$
that
\begin{align}
f_{m}(n)&=f_{m}(0)+(m+1)\sum_{i=1}^{n}f_{m-1}(n-i) \nonumber\\
&=f_{m}(0)+(m+1)\sum_{i=0}^{n-1}f_{m-1}(i). \label{recsquare}
\end{align}

We now give the natural generalization to Proposition
\ref{prop:attsheffpolcal}.
\begin{prop}
\label{prop:attsheffpolcalgen}
We have that $\mathcal{D}_{B}(m):\mathbb{C}[x] \rightarrow
\mathbb{C}[x]$, the generalized umbral shift attached to $B(t)$, is
characterized as the unique linear map satisfying:
\begin{align}
\label{eq:genpolcal}
\mathcal{D}_{B}(m)B_{n}(x)=f_{m}(n)B_{n-m}(x),
\end{align}
where by convention $B_{n}(x)=0$ for $n \leq -1$.
\end{prop}
\begin{proof}
By definition \ref{def:genattshef} and Lemma \ref{fsubmlemma} we have
\begin{align*}
\mathcal{D}_{B}(m)B_{n}(x)&=\phi_{B}L(m)L(-1)^{n}y\\
&=f_{m}(n)\phi_{B}L(-1)^{n-m}y\\
&=f_{m}(n)B_{n-m}(x).
\end{align*}
\end{proof}

\begin{remark} \rm
The convention in Proposition \ref{prop:attsheffpolcalgen} that
$B_{n}(x)=0$ for $n \leq -1$ is only used to ensure that in all cases
the right hand side of (\ref{eq:genpolcal}) is well defined.
This condition could have allowed for much more flexibility since we
already have that $f_{m}(n)=0$ whenever $n-m \leq -1$.
\end{remark}

We shall next solve for simple closed (polynomial) expressions for
$f_{m}(n)$ for fixed $m$.  If we are willing to sum over squares,
cubes etc. we could compute the answer for nonnegative $n$ for each
(fixed) $m$ in turn, using (\ref{recsquare}).  It is easy to verify
that $f_{0}(n)=n+1/2$ and $f_{1}(n)=n^{2}$ solve the first two cases.
To solve for the remaining cases however, for variety, we shall use a
heuristic argument making use of the Virasoro algebra relations to
derive the answer.  We have, for $m+l \leq n$,
\begin{align*}
(l-m)f_{l+m}(n)L(-1)^{n-l-m}y
&=(l-m)L(l+m)L(-1)^{n}y\\
&=[L(l),L(m)]L(-1)^{n}y\\
&=L(l)L(m)L(-1)^{n}y-L(m)L(l)L(-1)^{n}y\\
&=L(l)f_{m}(n)L(-1)^{n-m}y-L(m)f_{l}(n)L(-1)^{n-l}y\\
&=f_{l}(n-m)f_{m}(n)L(-1)^{n-m-l}y\\
&\quad -f_{m}(n-l)f_{l}(n)L(-1)^{n-l-m}y\\
&=\left(f_{l}(n-m)f_{m}(n)-f_{m}(n-l)f_{l}(n)\right)L(-1)^{n-l-m}y.
\end{align*}
We shall for the time being (unmathematically) ignore the restriction
on the indices and get, for whenever all terms are well defined, the
identity
\begin{align*}
(l-m)f_{l+m}(n)=f_{l}(n-m)f_{m}(n)-f_{m}(n-l)f_{l}(n).
\end{align*}
It is easy to see that the case $l=-1$ recovers $\ref{rec}$.  Further,
one can check that the case $l=0$ does not yield any new information.
For $l=1$ we get
\begin{align*}
(1-m)f_{m+1}(n)&=f_{1}(n-m)f_{m}(n)-f_{m}(n-1)f_{1}(n)\\
&=(n-m)^{2}f_{m}(n)-n^{2}f_{m}(n-1),
\end{align*}
so that it is easy to calculate, by simple substitution, each higher
case (in $m$) starting with $m=2$.  The calculations are not
difficult, of course, but I myself ``cheated'' and used Maple to find
and factor the first few answers, which yield an easy and obvious
pattern as follows:
\begin{align*}
f_{-1}(n)&=1\\
f_{0}(n)&=n+1/2\\
f_{1}(n)&=n^{2}\\
f_{2}(n)&=(1/2)n(n-1)(2n-1)\\
f_{3}(n)&=n(n-1)^{2}(n-2)\\
f_{4}(n)&=(1/2)n(n-1)(n-2)(n-3)(2n-3)\\
f_{5}(n)&=n(n-1)(n-2)^{2}(n-3)(n-4)\\
f_{6}(n)&=(1/2)n(n-1)(n-2)(n-3)(n-4)(n-5)(2n-5)\\
f_{7}(n)&=n(n-1)(n-2)(n-3)^{2}(n-4)(n-5)(n-6)\\
&\,\,\,\, \vdots
\end{align*}
With that as a guide, we may now return to doing rigorous math and
state and prove the following theorem.
\begin{theorem}
\label{the:fmn}
The unique solution to the recurrence equation (\ref{rec}) with $m
\geq -1$ and $n \in \mathbb{Z}$ and with boundary given by
$f_{-1}(n)=1$, $f_{0}(0)=1/2$ and $f_{m}(0)=0$ for $m \geq 1$ is given
by
\begin{align*}
f_{-1}(n)&=1\\
f_{m}(n)&=(1/2)(n(n-1)(n-2) \cdots
(n-m+1))(2n-m+1) \quad \text{ for} \qquad m \geq 0.
\end{align*}
\end{theorem}
\begin{proof}
The proof is a straightforward calculation.  Let $m \geq 0$ (although
admittedly the low $m$ cases are a bit degenerate in this notation).
Then we have
\begin{align*}
f_{m}(n)-f_{m}(n-1)&=(1/2)(n(n-1)(n-2) \cdots (n-m+1))(2n-m+1)\\
&\quad-(1/2)(n-1)(n-2) \cdots (n-m)(2(n-1)-m+1)\\
&=(1/2)(n-1)(n-2) \cdots (n-m+1)\cdot\\
&\quad\cdot(n(2n-m+1)-(n-m)(2(n-1)-m+1))\\
&=(1/2)(n-1)(n-2) \cdots (n-m+1)(2mn+2n-m^{2}-m)\\ 
&=(n-1)(n-2) \cdots (n-m+1)(m+1)(n-m/2)\\
&=(m+1)f_{m-1}(n-1).
\end{align*}
\end{proof}

\begin{remark} \rm
Recalling (\ref{recsquare}), it is easy to see that we could use the
last result, perhaps somewhat awkwardly, to solve for the sum of
squares and cubes etc., which happens to be related to the Bernoulli
numbers, one of the motivating subjects for Blissard \cite{Bli} and is
one of the classic problems solved via umbral methods (cf. Chapter 11
\cite{Do} for a nice, succinct old-fashioned umbral style proof and
also Chapter 3.11 \cite{Mel}).
\end{remark}

\begin{remark} \rm
We note that the umbral calculus has long been known to have
connections to the Bernoulli numbers and polynomials (see
e.g. \cite{Mel}).  Bernoulli polynomials have also appeared in the
literature of vertex algebra theory (see e.g. \cite{L1} and
\cite{DLM}).  Just as we have been establishing some analogues and
connections between umbral calculus and vertex algebra theory, it
might be interesting in future work to investigate further possible
connections explicitly related to Bernoulli numbers and polynomials.
\end{remark}

\begin{remark} \rm
We used a heuristic argument emphasizing the connection between umbral
calculus and the Virasoro algebra to help guess a solution to
$f_{m}(n)$ leading up to Theorem \ref{the:fmn}.  However, as the
referee has pointed out, one can connect the result with classical
umbral calculus as well.  The calculations below essentially reproduce
results obtained and shown to me by the referee and we shall follow
their reasoning.  Any errors or inelegance etc. in the particular
exposition presented here are entirely due to me.  Observing in all
that follows the typical caveats about the degeneracy of cases for
``low values of $n$'', for $n \geq 0$ let
\begin{align}
\label{eq:tf}
t_{n}(x)=f_{n-1}(x+n).
\end{align} 
Assuming we know that $t_{n}(x)$ are polynomials, then using
(\ref{rec}) for $n \geq 1$ we get
\begin{align*}
t_{n}(x)=t_{n}(x-1)+nt_{n-1}(x),
\end{align*}
since the result holds for all positive integral values of $x$.  Using
the formal Taylor theorem (\ref{translationoperator}), we get for $n
\geq 1$ that
\begin{align*}
\left(1-e^{-\frac{d}{dx}}\right)t_{n}(x)&=t_{n}(x)-t_{n}(x-1)\\
&=t_{n}(x-1)+nt_{n-1}(x)-t_{n-1}(x-1)\\
&=nt_{n-1}(x).
\end{align*}
Therefore $t_{n}(x)$ is a Sheffer sequence (cf. Theorem 2.3.7
\cite{Rm1}).  It is easy to see (cf. Theorem 2.4.3 in \cite{Rm1}) that
the relevant associated sequence $p_{k}(y)$ satisfies
\begin{align*}
\sum_{k \geq 0}\frac{p_{k}(y)}{k!}t^{k}=e^{y \log (1/(1-t))}=(1-t)^{-y},
\end{align*}
which yields 
\begin{align*}
p_{k}(y)=y(y+1)\cdots(y+k-1).
\end{align*}
We further have (cf. Theorem 2.3.9 \cite{Rm1}) 
\begin{align*}
t_{n}(x+y)=\sum_{k=0}^{n}\binom{n}{k}p_{k}(y)t_{n-k}(x).
\end{align*}
It is easy to see that $t_{0}(-2)=1$, $t_{1}(-2)=-1/2$ and that
$t_{n}(-2)=0$ for $n > 1$ so that for $n \geq 0$
\begin{align}
t_{n}(x)&=\sum_{k=0}^{n}\binom{n}{k}p_{k}(x+2)t_{n-k}(-2) \nonumber\\
&=p_{n}(x+2)-\frac{n}{2}p_{n-1}(x+2) \nonumber \\
&=(x+2)(x+3)\cdots(x+n+1)-\frac{n}{2}(x+2)(x+3)\cdots(x+n) \nonumber\\
&=\frac{1}{2}(x+2)(x+3)\cdots(x+n)(2x+n+2). \label{eq:tnx}
\end{align}
It is easy to check this formula using (\ref{eq:tf}) and Theorem
\ref{the:fmn}.  The referee also pointed out that letting
\begin{align*}
s_{n}(x)=\frac{t_{n}(x)}{n!}\left(=\frac{f_{n-1}(x+n)}{n!}\right),
\end{align*}
it is easy to see that the following particularly simple recursion holds
\begin{align*}
s_{n}(x)=s_{n-1}(x)+s_{n}(x-1),
\end{align*}
and, continuing to follow the referee, using (\ref{eq:tnx}) it is easy
to see that
\begin{align*}
s_{n}(x)=\binom{x+n+1}{n}-\frac{1}{2}\binom{x+n}{n-1}.
\end{align*}
\end{remark}

We shall conclude this paper by stating and proving the natural
generalization to the original formula defining the attached umbral
shifts in Definition \ref{def:umbshift}.
\begin{prop}
\label{prop:genumbshiftan}
For each $m \geq -1$, the map $\mathcal{D}_{B}(m):\mathbb{C}[x]
\rightarrow \mathbb{C}[x]$ is the unique linear map satisfying:
\begin{align*}
\mathcal{D}_{B}(m)e^{xB(w)}=
\left(w^{m+1}\frac{\partial}{\partial w}+\frac{m+1}{2}w^{m}\right)e^{xB(w)}.
\end{align*}
\end{prop}
\begin{proof}
We calculate, for $m \geq -1$ (although once again the low $m$ cases
are a bit degenerate) to get
\begin{align*}
w^{m+1}\frac{\partial}{\partial w}e^{xB(w)}
&=w^{m+1}\frac{\partial}{\partial w}\sum_{n \geq
0}\frac{B_{n}(x)w^{n}}{n!}\\ 
&=\sum_{n \geq 0}\frac{nB_{n}(x)w^{n+m}}{n!}\\
&=\sum_{n \geq m}\frac{(n-m)B_{n-m}(x)w^{n}}{(n-m)!}\\
&=\sum_{n \geq m}\frac{n(n-1)\cdots(n-m)B_{n-m}(x)w^{n}}{n!}
\end{align*}
and
\begin{align*}
w^{m}e^{xB(w)}&=\sum_{n \geq 0}\frac{B_{n}(x)w^{n+m}}{n!}\\ 
&=\sum_{n \geq m}\frac{B_{n-m}(x)w^{n}}{(n-m)!}\\
&=\sum_{n \geq m}\frac{n(n-1)\cdots(n-m+1)B_{n-m}(x)w^{n}}{n!},
\end{align*}
so that it is easy to check that
\begin{align*}
\left(w^{m+1}\frac{\partial}{\partial w}+\frac{m+1}{2}w^{m}\right)e^{xB(w)}&=
\sum_{n \geq m}\frac{f_{m}(n)B_{n-m}(x)w^{n}}{n!}\\
&=\mathcal{D}_{B}(m)e^{xB(w)}.
\end{align*}
\end{proof}

\begin{remark} \rm
\label{rem:Dan2}
Building on Remarks \ref{rem:Dan0} and \ref{rem:Dan1} we may regard
Proposition \ref{prop:genumbshiftan} as an analogue of formula
(8.7.37) in \cite{FLM} in the cases where the $n$ in (8.7.37) in
\cite{FLM} is restricted so that $n \geq -1$.  The $h$ in formula
(8.7.37) in \cite{FLM} should be replaced by the weight of the
relevant lowest weight vector, which in our setting seems perhaps to
correspond with the lowest weight of the module of the Virasoro
algebra which we have been considering, which as we have noted is
indeed $1/2$.
\end{remark}

\noindent {\small \sc Department of Mathematics, Rutgers University,
Piscataway, NJ 08854} 
\\ {\em E--mail
address}: thomasro@math.rutgers.edu

\end{document}